\documentclass[12pt]{article}
\usepackage[final]{epsfig}
\usepackage{graphics}
\usepackage{graphicx, amsfonts,amssymb, 
	amsmath,amsthm,color,mathtools,tikz,tikz-3dplot,hyperref,url}

\usepackage{latexsym}
\usepackage{amssymb}
\usepackage{graphicx}
\usepackage{url}
\usepackage{epstopdf}
\usepackage{wasysym}
\usepackage{color}
\usetikzlibrary{calc,intersections,through,backgrounds, arrows,decorations.markings,arrows.meta}
\usepackage{tkz-euclide}

\theoremstyle{definition}

\begin{document}

\def\mainfile{}
\newcommand{\eps}{{\varepsilon}}
\newcommand{\C}{{\mathbb C}}
\newcommand{\Q}{{\mathbb Q}}
\newcommand{\R}{{\mathbb R}}
\newcommand{\Z}{{\mathbb Z}}
\newcommand{\RP}{{\mathbb {RP}}}
\newcommand{\CP}{{\mathbb {CP}}}
\newcommand{\Tr}{\rm Tr}
\newcommand{\g}{\gamma}
\newcommand{\G}{\Gamma}
\newcommand{\e}{\varepsilon}

\title{Open problems on billiards and geometric optics}

\author{Misha Bialy\footnote{
School of Mathematical Sciences, 
Raymond and Beverly Sackler Faculty of Exact Sciences,
 Tel Aviv University, Israel; 
 bialy@tauex.tau.ac.il}
 \and
 Corentin Fierobe\footnote{
IST Austria, Am Campus 1, Klosterneuburg, Austria; 
corentin.fierobekoz@gmail.com}
\and
Alexey Glutsyuk\footnote{
CNRS (UMPA, ENS de Lyon), France, and HSE University, Moscow, Russia;
 aglutsyu@ens-lyon.fr
}
\and
Mark Levi\footnote{
Department of Mathematics,
Pennsylvania State University,
USA;
mxl48@psu.edu}
\and
Alexander Plakhov\footnote{
Department of Mathematics, 
University of Aveiro, Portugal;
plakhov@ua.pt}
\and
Serge Tabachnikov\footnote{
Department of Mathematics,
Pennsylvania State University,
USA;
sot2@psu.edu}
}

\date{\today}
\maketitle

This is a collection of problems composed by some participants of the workshop 
``Differential Geometry, Billiards, and Geometric Optics" that took place at CIRM on October 4--8, 2021.
\vskip 10 mm 

\centerline {\large\bf  Misha Bialy\footnote{Supported by ISF grant 580/20.}}

\begin{figure}[h]
	\centering
	\includegraphics[width=0.35\linewidth]{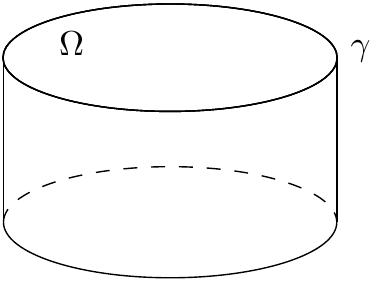}
	\caption{The surface.}
	\label{puck}
\end{figure}

\paragraph{Problem 1.}
Let $\gamma $ be a closed smooth  strictly convex curve in the plane bounding a domain $\Omega$, for example, an ellipse. Consider the geodesic flow on the surface homeomorphic to sphere, which is a cylinder of  height $d$ with base $\gamma$, glued on the top and on the bottom to $\Omega$. On the top and the bottom of the surface the motion is along straight lines, and on the cylindrical part along geodesics of the cylinder, see Figure \ref{puck}. 

We can describe this motion as follows.
Introduce arc-length coordinate $s$ on the curve $\gamma$ and let $\mathbf A$ be the phase cylinder with the coordinates $(s,\alpha)$. 

A line starting at $\gamma(s)$  with an angle $\alpha$ on the top comes to the boundary with an angle $\alpha_1$ at $\gamma(s_1)$. Then it passes to the surface of the cylinder with the same angle $\alpha_1$ and travels along the geodesic of the cylinder until it hits the bottom at $\gamma(s_1 +d\cot\alpha_1)$ with the same angle  $\alpha_1$. Next it passes to the bottom domain with the angle $\alpha_1$ and so on.

This geodesic motion can be described as a map $T$ of the cylinder $\mathbf A$ as follows.
$$
T:(s,\alpha)\mapsto (s_1+d\cot\alpha_1, \alpha_1).
$$
This means that $T$ is a composition of two maps, $T=T_2 \circ T_1,$ where $T_1$ is the usual billiard map $T_1:(s,\alpha)\mapsto (s_1, \alpha_1)$ and $T_2:(s,\alpha_1)\mapsto (s+d\cot\alpha_1 , \alpha_1)$ is the shift of $s$-coordinate.
Notice that $T$ is a symplectic map of the cylinder. Notice also that if $\alpha_1$ is small, then $d\cot\alpha_1$  is a large shift.
We ask

{\it Question 1}: Are there invariant curves of $T$?   For example, are there KAM curves near the boundary?

{\it Question 2}: What are the shapes of $\gamma$ (other than the circles) such that $T$ is an integrable map?

{\it Question 1}: Can $T$ be ergodic?

\paragraph{Problem 2.}

This is an old question about the integrability of outer billiards, see \cite{Ta} for further discussions and results. 

\begin{figure}[h]
	\centering
	\includegraphics[width=0.6\linewidth]{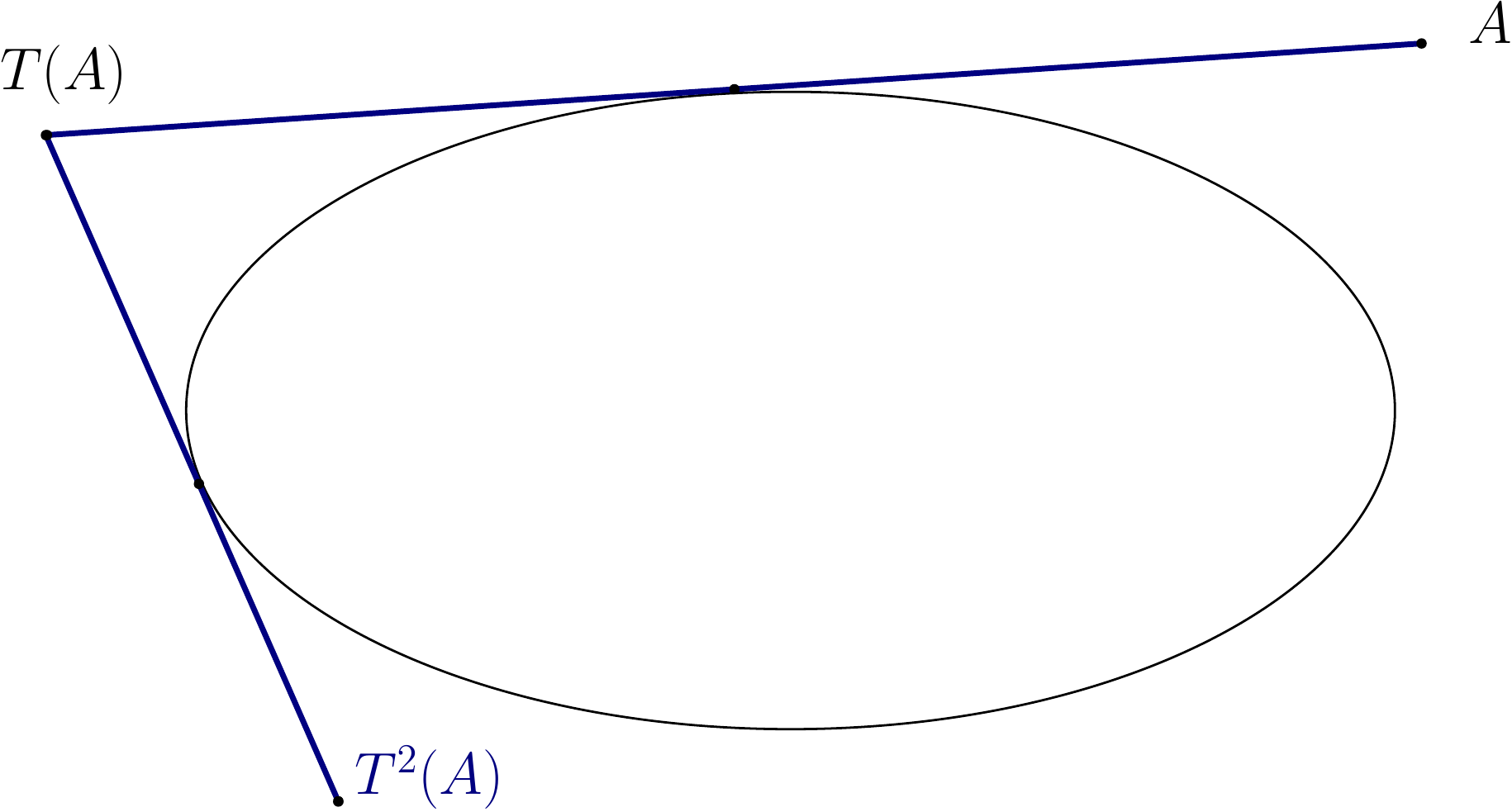}
	\caption{Outer billiard.}
	\label{outer}
\end{figure}

Let $\gamma$ be a smooth closed strictly convex curve in the plane. Let $T$ be the outer billiard map acting in the vicinity
of $\gamma$, see Figure \ref{outer}. The outer billiard reflection law reads:  the segment $[A,T(A)]$ is tangent to $\gamma$ at the midpoint. 

Obviously, outer billiard about an ellipse is integrable. Namely, in this case the whole phase space is foliated by the homothetic ellipses which are invariant curves for the outer billiard map.

{\it Question}: Are there other integrable outer billiards? 

This question can be considered as analogous to Birkhoff's conjecture (see recent paper 
\cite{BM} and the references therein)  for the usual billiards. 

\paragraph{Problem 3.}
This problem is about Gutkin billiards on the sphere and the hyperbolic plane. 

Let $\gamma$ be a smooth closed strictly convex curve in the plane, different from a circle, and let $\delta\in(0,\pi/2)$. We say that $\gamma$ has $\delta$-Gutkin property if the curve consisting of the incoming oriented lines  with the constant angle $\delta$  is an invariant curve of the Birkhoff billiard inside $\gamma$. 

E. Gutkin studied those billiards in the plane \cite{Gu}. Here we ask what  the Gutkin billiards are on the sphere and the hyperbolic plane. A first step was made in \cite{Tab}, where  infinitesimal deformations of the circle with the Gutkin property were studied. 

\paragraph{Problem 4.} Consider Birkhoff billiard inside a closed smooth strictly convex hypersurface $S\subset\mathbb{R}^d$. In the case $d=2$, Lazutkin proved  the existence of large family of convex caustics near $S$. However, for $d\geq 3$, caustics may exist only for ellipsoids \cite{berger, Be, gruber}. Nevertheless, it still can happen that the billiard map $T$, acting on the space $\mathbb A$ of all oriented lines intersecting $S$, leaves invariant a smooth hypersurface $\Sigma\subset\mathbb A$. 

It would be interesting to have an example of the billiard table (different from ellipsoids) $S\subset\mathbb{R}^d$ having an invariant hypersurface $\Sigma\subset\mathbb A$. What are geometric/dynamical properties of those $\Sigma$?

\paragraph{Problem 5.}  Consider a convex, smooth billiard table $\gamma$ in the plane which is symmetric with respect to a symmetry axis $l$. Suppose that $C$ is a convex caustic of the billiard. Prove, or give a counterexample, that $C$ is necessarily symmetric with respect to $l$.
One can show that if a non-symmetric caustic $C$ exists, then the rotation number of the corresponding invariant curve must be rational.
This problem appeared in \cite{arnold}.

\vskip 10 mm

\centerline {\large\bf  Corentin Fierobe}

\paragraph{Problem 1.} The following problem is related to a wide class of billiards, called \textit{projective billiards}. They were introduced and studied by Tabachnikov \cite{taba_projectif_ball, taba_projectif}. 

\begin{figure}[!ht]
\centering
\begin{tikzpicture}[line cap=round,line join=round,>=triangle 45,x=1.0cm,y=1.0cm]
\clip(-1.6,-1.2) rectangle (1.6,1.2);
\draw [rotate around={0:(0,0)}] (0,0) ellipse (1.41cm and 1cm);
\draw (-1.56,-0.1)-- (-1.27,0.11);
\draw (-1.51,0.3)-- (-1.12,0.39);
\draw (-1.24,0.73)-- (-0.96,0.53);
\draw (-0.88,0.99)-- (-0.64,0.69);
\draw (-0.3,1.18)-- (-0.2,0.79);
\draw (0.4,1.16)-- (0.31,0.78);
\draw (1.02,0.91)-- (0.77,0.64);
\draw (1.37,0.58)-- (1.09,0.42);
\draw (1.55,0.17)-- (1.28,-0.06);
\draw (1.07,-0.87)-- (0.76,-0.65);
\draw (0.66,-1.09)-- (0.37,-0.77);
\draw (0.14,-1.2)-- (0,-0.8);
\draw (-0.46,-1.15)-- (-0.39,-0.76);
\draw (-0.81,-1.03)-- (-0.79,-0.63);
\draw (-1.11,-0.85)-- (-1.03,-0.47);
\draw (-1.4,-0.54)-- (-1.24,-0.19);
\draw (1.11,-0.39)-- (1.52,-0.26);
\begin{scriptsize}
\end{scriptsize}
\end{tikzpicture}
\hspace{1cm}
\begin{tikzpicture}[line cap=round,line join=round,>=triangle 45,x=0.5cm,y=0.5cm]
\clip(-3.4,-2.5) rectangle (2.4,4);
\draw (2,1)-- (0,-2);
\draw (0,-2)-- (-3,3);
\draw (-3,3)-- (2,1);
\draw (-2.84,3.65)-- (-1.98,1.94);
\draw (-1.92,3.28)-- (-1.34,1.68);
\draw (-1.14,2.97)-- (-0.79,1.46);
\draw (-0.43,2.69)-- (-0.3,1.27);
\draw (0.2,2.43)-- (0.14,1.09);
\draw (0.86,2.17)-- (0.6,0.9);
\draw (1.58,1.88)-- (1.1,0.7);
\draw (-2.99,1.75)-- (-1.37,1.51);
\draw (-2.19,0.42)-- (-0.83,0.61);
\draw (-1.43,-0.84)-- (-0.32,-0.25);
\draw (-0.79,-1.91)-- (0.12,-0.97);
\draw (0.07,-0.71)-- (0.96,-1.79);
\draw (1.36,-1.2)-- (0.38,-0.25);
\draw (0.68,0.21)-- (1.75,-0.6);
\draw (0.96,0.62)-- (2.1,-0.07);
\draw [line width=0.4pt,dash pattern=on 1pt off 1pt] (0.97,-0.82)-- (-3,3);
\draw [line width=0.4pt,dash pattern=on 1pt off 1pt] (-1.22,0.55)-- (2,1);
\draw [line width=0.4pt,dash pattern=on 1pt off 1pt] (0.15,1.3)-- (0,-2);
\end{tikzpicture}
\hspace{1cm}
\begin{tikzpicture}[line cap=round,line join=round,>=triangle 45,x=0.5cm,y=0.5cm]
\clip(-3.7,-2.5) rectangle (3.7,3.5);
\draw (-2,-2)-- (-3,1);
\draw (-3,1)-- (2,3);
\draw (2,3)-- (3,-2);
\draw (-2,-2)-- (3,-2);
\draw [dotted] (-2,-2)-- (2,3);
\draw [dotted] (-3,1)-- (3,-2);
\draw [dash pattern=on 1pt off 1pt] (-0.57,-0.21)-- (-2.94,-0.84);
\draw [dash pattern=on 1pt off 1pt] (-0.57,-0.21)-- (-0.48,2.46);
\draw [dash pattern=on 1pt off 1pt] (-0.57,-0.21)-- (2.95,0.71);
\draw [dash pattern=on 1pt off 1pt] (-0.57,-0.21)-- (0.4,-2.46);
\draw (-3.39,0.53)-- (-2.2,0.21);
\draw (-3.17,-0.15)-- (-2.06,-0.18);
\draw (-2.94,-0.84)-- (-1.93,-0.58);
\draw (-2.68,-1.62)-- (-1.78,-1.02);
\draw (-1.5,-2.46)-- (-1.11,-1.51);
\draw (-0.57,-1.51)-- (-0.56,-2.46);
\draw (-0.01,-1.51)-- (0.4,-2.46);
\draw (1.25,-2.46)-- (0.48,-1.51);
\draw (1.06,-1.51)-- (2.25,-2.46);
\draw (3.33,-1.23)-- (2.3,-0.96);
\draw (3.2,-0.55)-- (2.2,-0.46);
\draw (3.07,0.1)-- (2.1,0.01);
\draw (2.01,0.46)-- (2.95,0.71);
\draw (2.82,1.34)-- (1.92,0.93);
\draw (1.82,1.44)-- (2.68,2.04);
\draw (2.55,2.68)-- (1.72,1.92);
\draw (1.36,3.19)-- (0.7,2.03);
\draw (0.53,2.86)-- (0.15,1.81);
\draw (-0.48,2.46)-- (-0.51,1.54);
\draw (-0.96,1.37)-- (-1.16,2.19);
\draw (-1.96,1.87)-- (-1.48,1.16);
\draw (-1.95,0.97)-- (-2.65,1.58);
\end{tikzpicture}
\caption{Left: a convex closed curve with a field of transverse lines. Center: the so-called right-spherical billiard which is $3$-reflective. It is defined by triangular domain, and each transverse line on a side join the opposite vertex of the triangle. Right: a $4$-reflective projective billiard inside a quadrilateral. The transverse lines join the intersection point of the latter's diagonals.}
\label{fig1}
\end{figure}
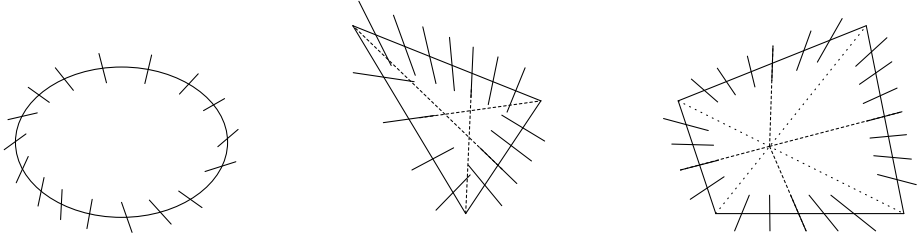

A projective billiard in the Euclidean plane can be defined as a bounded domain with a (piecewise) smooth boundary which is endowed with a smooth transverse line field (see Figure \ref{fig1}). Given such a domain $\Omega\subset\R^2$, whose boundary $\partial\Omega$ is endowed with a transverse line field $L$, we  define a \textit{law of reflection} of oriented lines intersecting $\partial\Omega$ as follows: if $\ell$ is an oriented line intersecting $\partial\Omega$ at a point $p$, we say that it is reflected into a line $\ell'$ intersecting $\partial\Omega$ at $p$ if the four lines $\ell$, $\ell'$, $L(p)$ and $T_p\partial\Omega$ form a \textit{harmonic} quadruple, that is, if the unique non-trivial involution of the pencil of lines containing $p$ fixing $L(p)$ and $T_p\partial\Omega$ permutes $\ell$ and $\ell'$.

Moreover, the line $\ell'$ is naturally oriented in the opposite direction to $\ell$. Then, by construction, when $L(p)$ is orthogonal to $\partial\Omega$, $\ell'$ is given by the usual reflection law of optics (angle of incidence $=$ angle of reflection), and therefore the class of projective billiards contains all usual billiards. 

\textbf{$k$-reflective billiards.} A projective billiard is said to be \textit{$k$-reflective} if its corresponding billiard map has an open set of periodic points. Unlike the usual billiards, it is not difficult to find examples of $k$-reflective projective billiards, namely for $k=3$ or for any even $k\geq 4$ (see Figure \ref{fig1}, and \cite{fierobe1, fierobe_these} for more precise descriptions and elementary proofs). These examples consist of projective billiards inside polygons. To our knowledge there are no examples of $k$-reflective projective billiard with $k\geq5$ odd. Therefore we ask the following questions.

{\it Question 1}: Can one find examples of $k$-reflective projective billiards with $k\geq 5$ odd?

{\it Question 2} (Projective Ivrii's conjecture): Given an integer $k\geq 4$, are there other examples of $k$-reflective projective billiards within a class of domains with a specific boundary smoothness (polygonal, piecewise-algebraic or analytic, etc.)?

\vskip 10 mm

\vskip 10 mm

\centerline {\large\bf  Alexey Glutsyuk\footnote{
Partially supported by Laboratory of Dynamical Systems and Applications, HSE University, of the Ministry of science and higher education of the RF grant ag. No 075-15-2019-1931 and by RFBR and JSPS (research project 19-51-50005).
}}

\paragraph{Problem 1.} The Birkhoff Conjecture concerns a bounded strictly convex planar billiard 
with smooth boundary. 

Suppose that some neighborhood of the boundary from the inner (convex) side 
 is foliated by closed caustics. The Birkhoff Conjecture states that in this case the billiard boundary is 
 an ellipse. 
 
\begin{figure}[ht]
  \centering
	\includegraphics[width=0.5\linewidth]{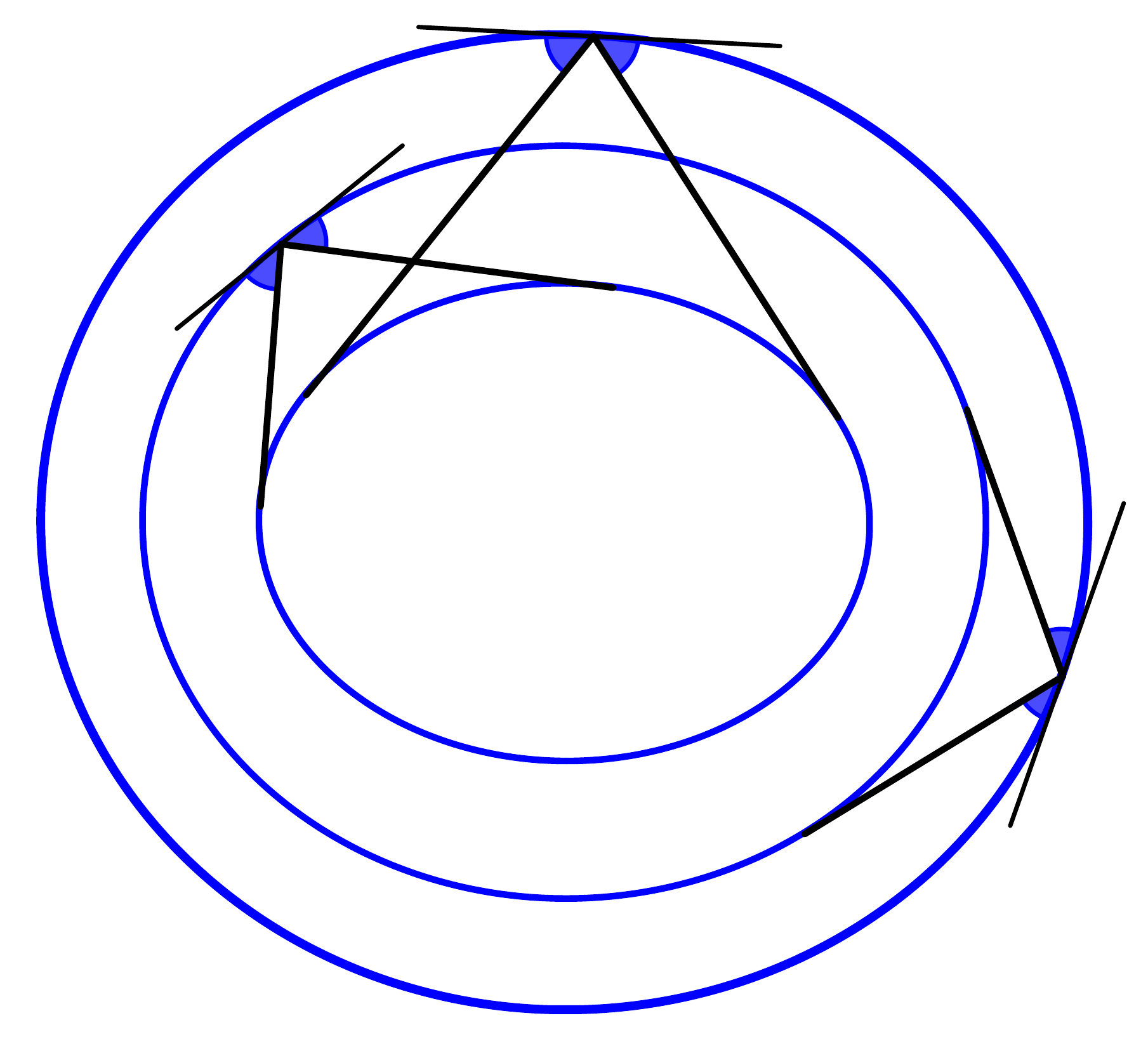}
	\caption{A caustic of a caustic is a caustic.}
	\label{curves}
\end{figure}

 One of the first famous results on this conjecture is due to Hillel Poritsky. 
 In his paper \cite{poritsky}, where he first stated the Birkhoff Conjecture in print, he proved it under 
 the additional assumption that for any two nested caustics the smaller caustic is a caustic for the 
 billiard in the bigger caustic. See the Figure \ref{curves}. 
 
 For  recent results and  survey of the Birkhoff Conjecture, see \cite{bm, kalsor, kalsor2}. 
 
 {\it Question} (Konstantin Khanin, August 2019): {Suppose that a bounded strictly convex planar billiard 
 has just {\it two} nested closed caustics satisfying the above condition:  the smaller caustic is a caustic for the billiard in the  bigger one. Is it true that the billiard boundary is an ellipse?}

\paragraph{Problem 2.} Consider a closed strictly convex hypersurface $\gamma\subset\mathbb R^n$. Let $\Pi$ denote 
 the phase cylinder for the billiard inside the domain bounded by $\gamma$: this is the set of those oriented 
 lines in $\mathbb R^n$ that intersect $\gamma$ transversally at two points.  It is an open subset in the space of 
 oriented lines equipped with the standard symplectic structure. 
 
 There  is an important open question stated in \cite{ptt}: 
which symplectomorphisms can be realized by compositions of reflections?
 
 {\it Question \cite[subsection 1.6]{gl}}: {Is it true that for every given $\varepsilon>0$ and  $k\in\mathbb N$ 
 each $C^{\infty}$-smooth Hamiltonian symplectomorphism $\Pi\to\Pi$ is the $C^{\infty}$-limit 
 of compositions of reflections from the hypersurface $\gamma$ and from hypersurfaces $\varepsilon$-close to 
 $\gamma$ in the $C^k$-topology?}
 \medskip
 
 {\it Remark.} Compositional differences of reflections from a hypersurface $\gamma$ and from 
 its deformations were  
 introduced in \cite{perline} and studied in \cite{perline, gl}. The positive answer to the above question was proved 
 in \cite{gl} for compositions of not  just reflections, but of reflections and of their inverses. 
 \medskip
 
The famous Ivrii's Conjecture \cite{Ivrii} states that in every 
 billiard with smooth boundary the set of periodic initial conditions for the billiard map has Lebesgue 
 measure zero. In particular, it implies a slightly weaker conjecture stating 
 that no iterate of the billiard ball map  can coincide with the identity on 
 an open subset in the space of oriented lines. 
 
 A bit stronger version of the latter conjecture states that 
 {\it no well-defined composition of reflections from smooth germs of  hypersurfaces can be equal to the identity.} 
 \medskip
 
{\it Remark.} 
If in this conjecture we replace ``composition of reflections" by ``composition of reflections and of inverses of 
 reflections", then the conjecture is false, since, for example, the billiard ball maps in confocal ellipses commute 
\cite[p.59, corollary 4.6]{Ta}, \cite[p.58]{tabcom}. 
\medskip

Due to this remark, it is important to understand which symplectomorphisms 
of the phase cylinder are limits of compositions of reflections, without including their inverses. 

\vskip 10 mm

\centerline {\large\bf  Mark Levi}

\paragraph{Problem 1.} Consider a curve $C$ in $ {\mathbb R}  ^2 $, and an $\varepsilon$-sized square lattice. As I translate the curve, every time it meets
a lattice point, a click sounds. Think of a line (translation parameter $\lambda$)  with points on it (the values of $\lambda$ when the curve contains a lattice point). 

{\it Question 1}: Can one recover the shape of the curve given this distribution of points for any $ \varepsilon>0 $ (presumably only small $\varepsilon$ are needed, of course) for the curve $C$? Or -- an easier question -- for the curve  and for all of its rotations? 

This question of recovering the shape of the curve from the cloud of clicks is a kind of a ``big data" question. It involves extracting/recognizing patterns in the cloud of indistinguishable clicks. These patters encode information about the curvature: for example, near points with ``very rational" slopes one can
write asymptotics of the distribution in terms of the curvature and the slope at the point of rational slope.  

{\it Question 2}: The problem is that all these clicks are dumped on the segment; is there a ``Fourier transform" which  extracts the ``curvature-induced patterns from a cloud of clicks? 

{\it Question 3}: Analyze special cases of the segment; a polygon; a circle. 
 (A trivial observation: the cloud
evolves as a periodic function of the translation parameter $\lambda$, of period $\varepsilon$, so only $ \lambda \in [0, \varepsilon ] $ needs to be considered). 

\vskip 10 mm

\centerline {\large\bf  Alexander Plakhov\footnote{Supported by CIDMA through FCT (Funda\c{c}\~ao para a Ci\^encia e a Tecnologia), ref. UIDB/04106/2020.}}

\paragraph{Problem 1: dimension of a trapped set.} 

Consider a smooth flow or cascade on $\mathbb R^2$ preserving the Lebesgue measure,
$$
g^t : \mathbb R^2 \to \mathbb R^2, \qquad g^{t+s} = g^t \circ g^s, \quad g^0 = \text{id}, \quad t \in \mathbb R \ \ \text{or} \ \ t \in \mathbb Z.
$$
A point $x \in \mathbb R^2$ is called {\it trapped}, if the positive semi-orbit $\{ g^t x, \ t \ge 0 \}$ is bounded and the negative semi-orbit $\{ g^t x, \ t \le 0 \}$ is unbounded. Let $T_g \subset \mathbb R^2$ be the set of trapped points of $g$. The Lebesgue measure of $T_g$ is zero.

{\it Question:} What is the maximum Hausdorff dimension of $T_g$? 

{\it Example.} Consider the flow $g$ induced by the differential equations $\dot x = -x$, $\dot y = y$. Here the dimension of $T_g = \mathbb R \times \{0\}$ is 1.

{\it Hint.} The answer depends on the smoothness of $g$. In particular, if no smoothness assumptions are made, the answer is as follows: the dimension of $T_g$ for a certain measurable cascade $g$ is 2 (the maximum possible value). In general, the answer should look like this: let $g$ be $C^m$ in $x$ (or in both $x$ and $t$), then the maximum dimension of $T_g$ is $d(m)$ (a value between 1 and 2 to be determined).

This problem arose from a similar problem for trapping sets in billiards.

\paragraph{Problem 2: motion in Newtonian aerodynamics.}

Consider a solid body in $\mathbb R^2$, that is, a compact domain $B \subset \mathbb R^2$ with piecewise smooth boundary, with the mass uniformly distributed in $B$. The total mass of $B$ is $m > 0$. Consider a uniform medium in $\mathbb R^2 \setminus B$ composed of point particles at rest. (There are infinitely many particles of zero mass; the total mass of particles contained in any domain $D \subset \mathbb R^2 \setminus B$ equals the area of $D$).

At the instant $t=0$ the body starts the (translational and/or rotational) motion. When colliding with the body's boundary, the particles are reflected elastically.\footnote{More precisely, let a particle hit the body at a regular point $x$ of its boundary. One needs to take an inertial reference system in which the point $x$ is at rest; in this system the reflection occurs according to the familiar billiard law.}

{\it Problem:} Suppose that $B$ has a simple shape: an ellipse, a triangle, a rod $B = [-1,\, 1] \times \{0\}$, etc. (The trivial case when the body is a disc is excluded.) Describe the motion of $B$ for $t \ge 0$.

{\it Related questions:} The rod $B = [-1,\, 1] \times \{0\}$ starts the rotational motion around its center (no translation). Will the total number of turns be finite or infinite? In the latter case, will the velocity of rotation go to zero? What is the asymptotic behavior of the velocity as $t \to +\infty$?

{\it More questions:} A centrally symmetric body (for example, an ellipse, a square, etc.) starts the rotational motion around its center (no translation). The same questions as in the previous paragraph.

\paragraph{Problem 3: equations of motion in a rarefied medium.}
({This problem is closely connected with the previous one}.)

A body moves freely in the framework of Newtonian aerodynamics in $\mathbb R^n$, $n \ge 1$ in a rarefied medium. 

{\it Problem:} Determine the equation(s) of dynamics and prove a/the theorem of existence and uniqueness.

{\it Particular case: $n=1$.} A massive point with mass 1 moves on the line in a medium composed of identical infinitesimal particles. A motion should be defined by the triple $\mu_t,\, X(t),\, P(t)$, where $\mu_t$ is a measure on $\mathbb R^2 \ni (x,v)$ describing the distribution of particles at the instant $t$, and $X(t)$ and $P(t) = X'(t)$ are the position and the velocity of the massive particle at the instant $t$.
\vskip 10 mm

\centerline {\large\bf  Serge Tabachnikov\footnote{Supported by NSF grant DMS-2005444.}}

\begin{figure}[h]
	\centering
	\includegraphics[width=0.7\linewidth]{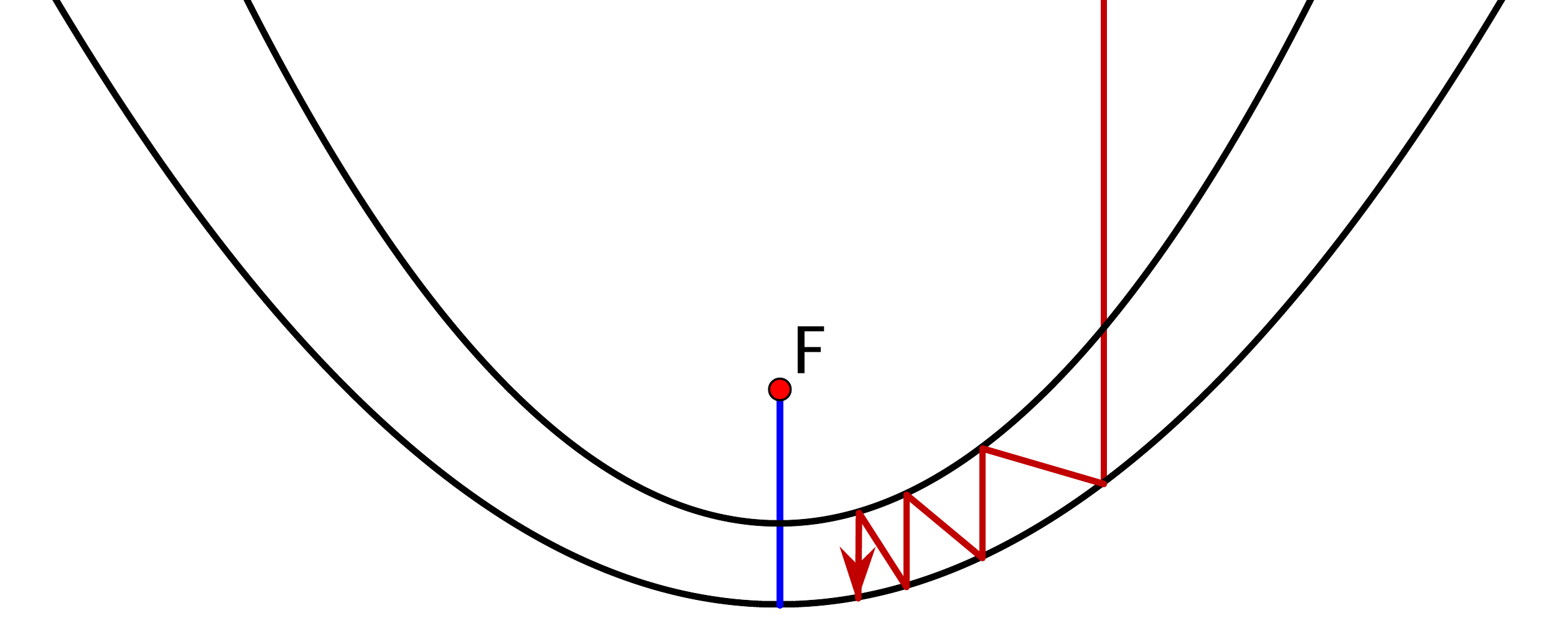}
	\caption{Trap for a parallel beam of light: $F$ is the focus of the parabolas.}
	\label{trap}
\end{figure}

\paragraph{Problem 1.} Figure \ref{trap} presents a trap for a planar vertical beam of light. Consider two confocal and coaxial parabolas, and let a ray of light, parallel to the axis, enter the domain between the parabolas through an aperture. The optical property of parabolas imply that if this ray optically reflects from the parabolas, it will be trapped inside the domain and will never cross the axis. See \cite{Ta} for another design of a trap.

The trap in  Figure \ref{trap} can capture a parallel beam of light. The space of oriented lines ${\mathcal L}$ in $\mathbb{R}^2$ is 2-dimensional, and a parallel beam is a curve in ${\mathcal L}$. Similarly, one can design a trap for any local 1-parameter family of rays (by first sending it to a parallel beam by reflection in a mirror). However, the Poincar\'e recurrence theorem implies that one cannot trap a 2-parameter family of rays.

Consider the situation in $\mathbb{R}^n$. Then dim ${\mathcal L}=2n-2$, and this space carries a natural symplectic structure. A variation of the previous construction yields a trap for a parallel beam and, more generally, for any local normal family of rays (the rays that are perpendicular to a hypersurface; such families comprise Lagrangian submanifolds in ${\mathcal L}$). 

{\it Question 1}: What is the greatest dimension of a family of rays in $\mathbb{R}^n$ that can be trapped?

{\it Question 2}: In $\mathbb{R}^3$, can one trap a non-normal 2-parameter family of rays?

\begin{figure}[h]
	\centering
	\includegraphics[width=0.4\linewidth]{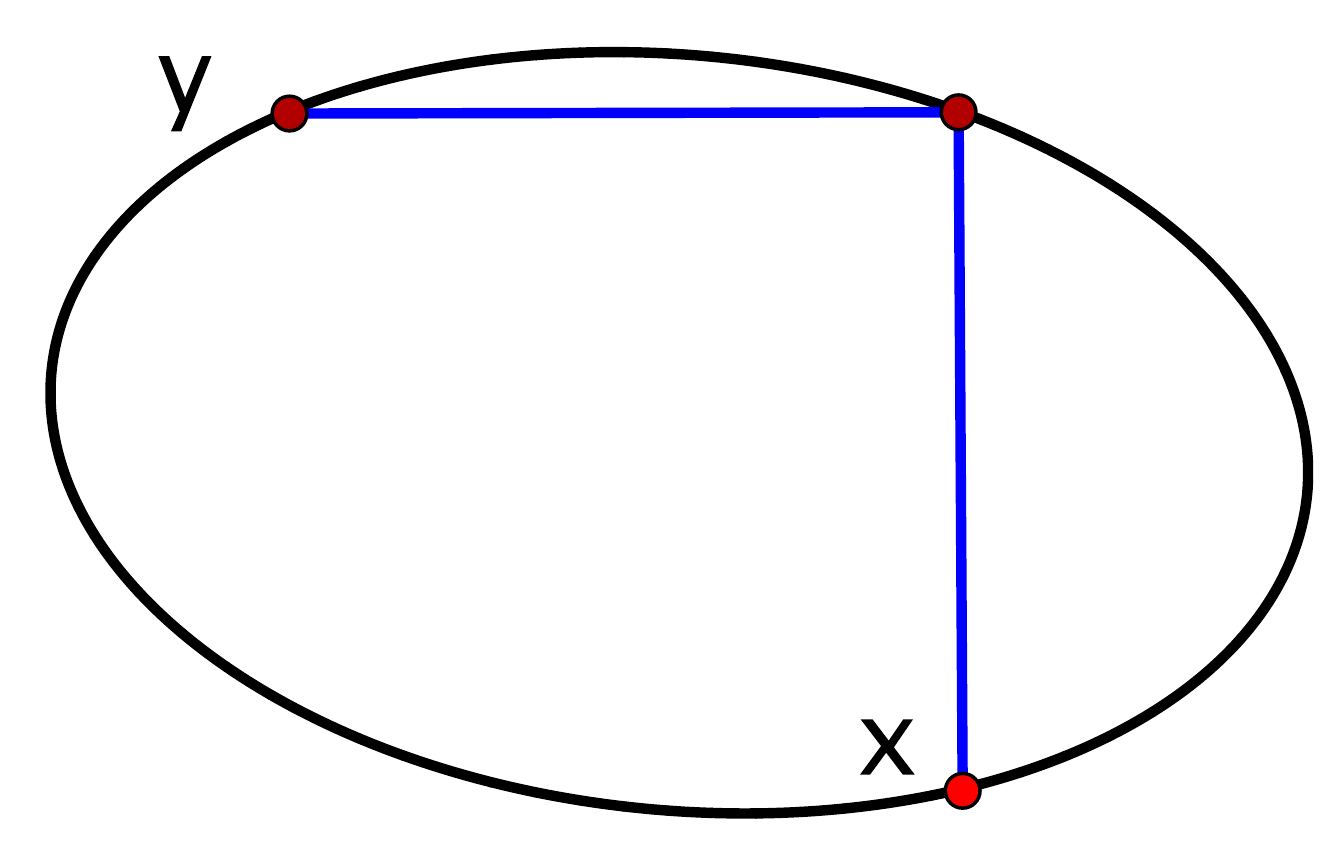}
	\caption{A circle map.}
	\label{oval}
\end{figure}

\paragraph{Problem 2.} Consider a planar oval $\gamma$ (a smooth closed strictly convex curve) and choose two directions, say, vertical and horizontal. We construct a map $F: \gamma\to \gamma$ as follows. Given a point $x\in\gamma$, draw the  line in the first direction through it until the intersection with $\gamma$, and then draw the line in the second direction through this new point until the intersection with $\gamma$ at point $y$. One has $F(x)=y$, see Figure \ref{oval}.

This map was considered in a number of papers from  different perspectives \cite{Ar,GKT,HM,Jo,KT,NT,So}.

If $\gamma$ is an ellipse, then $F$ is conjugated to a rotation for every choice of the pair of directions.

{\it Question 1}: Is this property characteristic of ellipses?

The answer is in the affirmative if one additionally assumes that $\gamma$ is centrally symmetric, see \cite{Ta21}.

Let us modify the construction by replacing the two families of parallel lines by two pencils of lines passing through points $P$ and $Q$. This provides a projective version of the previous map (where $P$ and $Q$ were points at infinity). 

If $\gamma$ is an ellipse and the line $PQ$ is disjoint from it, the respective map $F$ is still conjugated to a rotation. If 
an oval $\gamma$ has the property that for every choice of points $P$ and $Q$, such that $(PQ)\cap \gamma =\emptyset$, the map $F$ is conjugated to a rotation, then $\gamma$ is an ellipse, see \cite{Ta21}.

However if the line $PQ$ intersects $\gamma$, the map $F$ has two fixed points. If $\gamma$ is an ellipse, this map is conjugated to a M\"obius (projective) transformation.

{\it Question 2}: Let an oval $\gamma$ have the property that for every choice of points $P$ and $Q$, such that $(PQ)\cap \gamma \neq \emptyset$, the map $F$ is conjugated to a M\"obius transformation. Is it true that $\gamma$ is an ellipse?

\begin{figure}[h]
	\centering
	\includegraphics[width=0.35\linewidth]{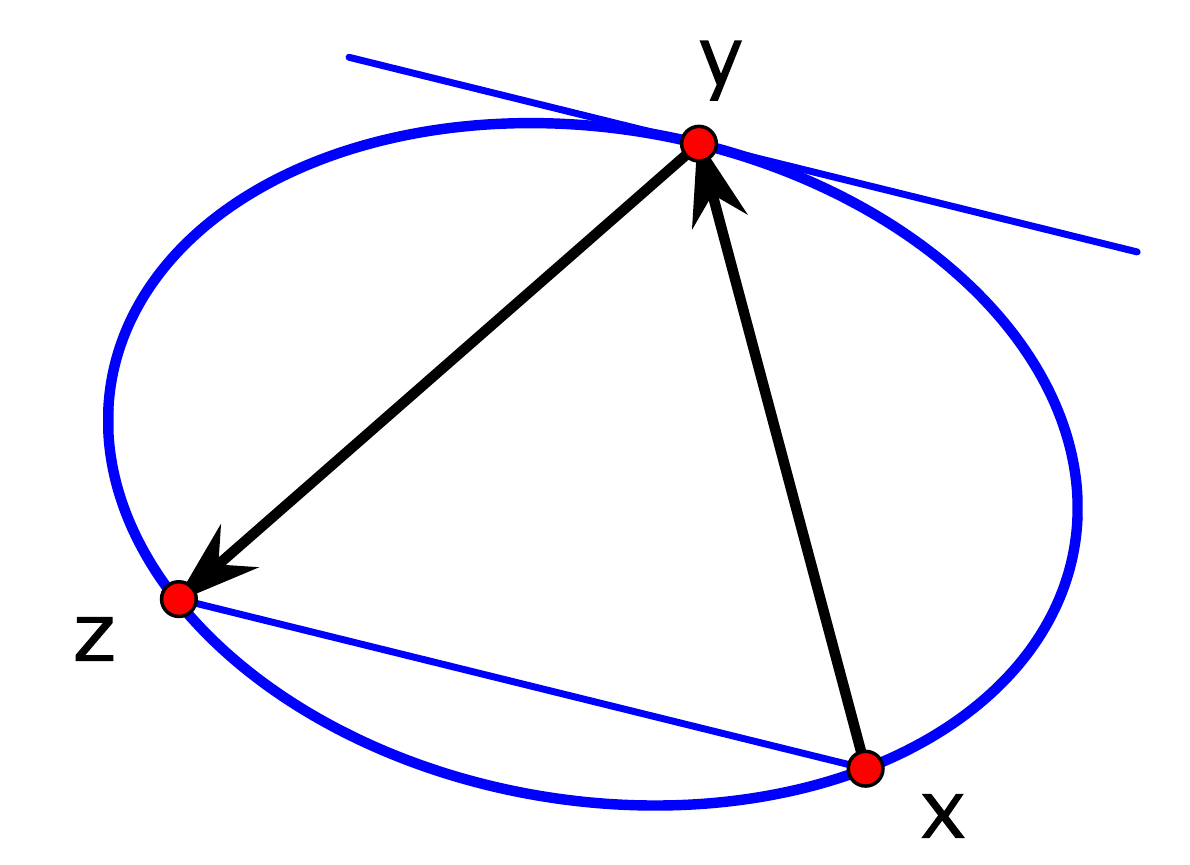}\qquad\quad
	\includegraphics[width=0.35\linewidth]{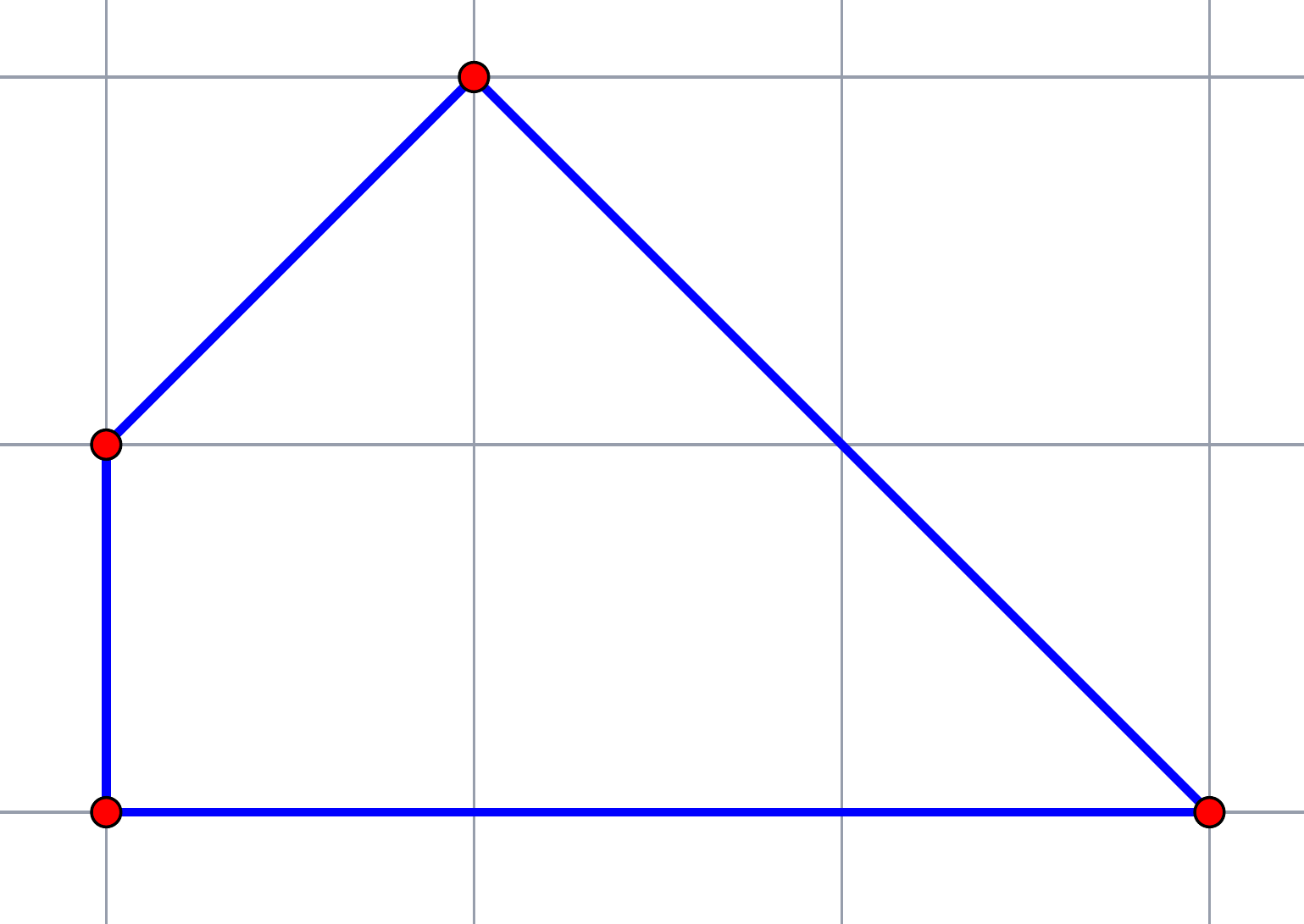}
	\caption{Left: the chord $xy$ reflects to $yz$ if the tangent line to the oval at $y$ is parallel to $xz$. Right: all orbits in this quadrilateral are periodic (with the periods equal to 20 and 36).}
	\label{reflection}
\end{figure}

\paragraph{Problem 3.} A planar symplectic billiard is a dynamical system on oriented chords of an oval given by a non-conventional reflection law depicted in Figure \ref{reflection}, see \cite{AT}. Similarly one defines polygonal symplectic billiards. 

In \cite{AT,ABSST}, a number of polygons are described that have the property that all symplectic billiard orbits are periodic (in particular, the affine-regular polygons and the trapezoids have this property).

{\it Question 1}: Describe all such polygons.

{\it Question 2}: Does every polygon have a periodic orbit?

For the usual billiards, the latter is a famous  problem, open even for triangles.

{\it Question 3}: Is the symplectic billiard dynamic in the stadium chaotic?

Numerically, this seems to be the case.

\begin{figure}[h]
	\centering
	\includegraphics[width=0.6\linewidth]{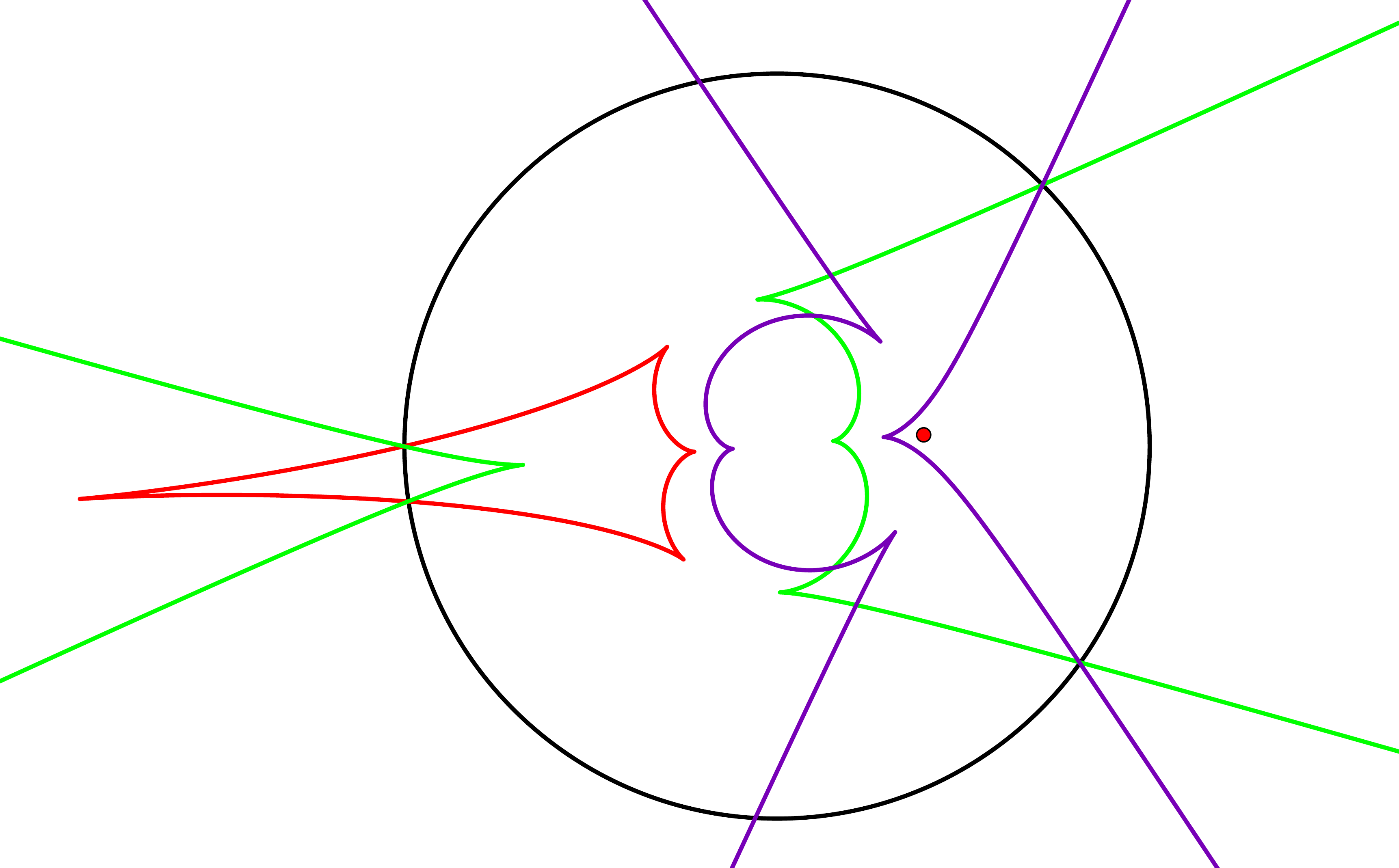}
	\caption{The first three caustics by reflection in a circle.}
	\label{three}
\end{figure}

\paragraph{Problem 4.} Consider an oval $\gamma$, thought of as an ideal mirror, and a source of light inside it. The envelope of the rays of light that have undergone $n$ reflections in $\gamma$ is called $n$th caustic by reflection. See Figure \ref{three}. 

One can prove that, for every oval, and the source of light in general position,  all caustics by reflection has at least four cusps, see \cite{BT}. 

{\it Question 1}: Is it true that each (generic) caustic by reflection in an ellipse has exactly four cusps?

{\it Question 2}: Does the above property characterize ellipses?

These questions are related to the Last Geometric Statement of Jacobi concerning the conjugate locus of a point on an ellipsoid, see, e.g., \cite{Si}.

\paragraph{Problem 5.}  The length spectrum of the billiard in an oval $\gamma$, that is, the set of lengths of its periodic trajectories, is related to the spectrum of the Laplace operator with the Dirichlet boundary condition in $\gamma$ by the Poisson summation formula, see \cite{GM}. 

The outer billiard about an oval $\gamma$ has the associated area spectrum: the set of areas of polygons formed by the periodic outer billiard trajectories (these are the circumscribed polygons of extremal areas). Is this area spectrum related to the spectrum of some differential operator?

Note that, in the spherical geometry, the inner and outer billiards are related by the spherical duality that assigns to an oriented great circle its pole  (see, e.g., \cite{Ta}).

\end{document}